\input amstex
\documentstyle{amsppt}
\pageheight{50.5pc} \pagewidth{32pc}

\define\1{\hbox{\rm 1}\!\! @,@,@,@,@,\hbox{\rm I}}

\topmatter
\title
{The small deviations of many-dimensional diffusion processes and
rarefaction by boundaries}
\endtitle

\author
{Vitalii A. Gasanenko}
\endauthor
\address
{Institute of Mathematics, National Academy of Science of Ukraine,
  Tereshchenkivska 3, 252601, Kiev, Ukraine}
\endaddress
\keywords{parabolic problem,small domain, algorithm of expansion,
number of unabsorbed processes}
\endkeywords
\subjclass {60 J 65}
\endsubjclass
\email {gs\@imath.kiev.ua or gsn\@ckc.com.ua}
\endemail
\abstract We lead the algorithm of expansion of sojourn probability
of many-dimensional diffusion processes in small domain. The
principal member of this expansion defines normalizing coefficient
for special limit theorems.
\endabstract
\endtopmatter
\rightheadtext{The small deviation of many-dimensional diffusion
processes}
\document

{\bf Introduction.}

\bigskip
Let $\xi(t)$ be a random process with measurable phase space
$(X,\Sigma(X))$. Consider the measurable connected domain $D\in
\Sigma(X)$ and small parameter $\epsilon$. The investigations of
asymptotics of sojourn probability (small deviations)

$$
P\left(\xi(t)\in \epsilon D,\quad t\in [0,T]\right)\eqno(1)
$$

is jointed with many practice and theoretical problems [1-4]. In the
literature, it was researched both rough asymptotics of  principal
member of (1)(log from it)[5] and exact asymptotics of diffusion
processes of (1)[6-8].  In the works [9,10] was proved of algorithms
of expansions of exact asymptotics of small deviation for diffusion
and piecewise deterministic random processes for one-dimensional
case.

The purpose this article is to present the algorithm of expansion of
small deviation for many-dimensional diffusion processes and to
define all constants of principal member.

In Section 1 our main result is stated and proved. In section 2 we
consider the limits theorems about numbers of unabsorbed diffusion
particles by boundaries of small domain.

\bigskip
{\bf I.} {\bf The expansion.}

\bigskip

We shall investigate of asymptote  of following probability

$$
P(\epsilon,x)=P\left(\xi(t)\in \epsilon D,\quad 0\leq t\leq
T\right),\quad \epsilon\to 0,
$$

where $\xi(t)\in R^{d}$ is solution of the following stochastic
differential equation

$$
d\xi(t)=a(t,\xi(t))dt + \sum\limits_{i=1}^{d}b_{i}(\xi(t))
dw_{i}(t), \quad \xi(0)= x\in \epsilon D.\eqno(2)
$$

where functions

$$
 b_{i}(x),~ a(t,x): R^{d}\to R^{d}\quad\hbox{and}\quad R_{+}\times R^{d}\to R^{d}.
$$

are differentiable.

Set $ \sigma_{ij}(x)=\sum\limits_{k}b_{k}^{i}(x)b_{k}^{j}(x)$.

It is known that $P(\epsilon,x)= u_{0}^{\epsilon}(T,x)$. Here
$u_{0}^{\epsilon}(t,x)$ is  solution of the following parabolic
boundary problem at $0\leq t\leq T$

$$
\frac{\partial u_{0}^{\epsilon}(t,x)}{\partial t}=
\frac{1}{2}\sum\limits_{i,j=1}^{d}\sigma_{ij}(x)\frac{\partial^{2}
u_{0}^{\epsilon}(t,x)} {\partial x_{i}\partial x_{j}}+
\sum\limits_{i=1}^{d}a_{i}(T-t,x) \frac{\partial
u_{0}^{\epsilon}(t,x)}{\partial x_{i}},
 \quad x\in  D_{\epsilon};
$$

$$
u(t,x)|_{t=0}=1;\quad x\in  D_{\epsilon};\quad u(t,x)=0\quad x\in
\partial D_{\epsilon},\quad 0\leq t\leq T.\eqno (3)
$$

where $D_{\epsilon}=\epsilon D$. It is assumed that  $D$ is a
connected bounded domain from $R^{m}$; the boundary $\partial Q$ is
the Lyapunov surface ~$C^{(1,\lambda)}$ and $0\in D$.  We interest
of the asymptotic expansion $\epsilon\to 0$ of solution this problem
$u_{0}^{\epsilon}(t,x)$ at $\epsilon\to 0$.

We define the differential operator $A:\frac{1}{2}\sum\limits_{1\leq
i,j\leq d}\sigma_{ij}(0)\frac{\partial^{2}} {\partial x_{i}\partial
x_{j}}.$ Let $\sigma$ be a matrix with the following property

$$
\sum\limits_{1\leq i,j\leq d}\sigma_{ij}(0)z_{i}z_{j} \geq \mu |\vec
z|^{2}.
$$

Here ~$\mu$,~ there is a fixed positive number, and $\vec
z=(z_{1},\cdots, z_{d})$~ is an arbitrary real vector.

This operator acts in the following space

$$
H_{A}=\{u: u\in L_{2}(D)\cap Au\in L_{2}(D)\cap u(\partial D)=0\}
$$
with inner product $(u,v)_{A}=(Au,v)$. Here $(,)$ is inner product
in $L_{2}(Q)$.
 The operator ~$A$~ is a positive
operator[11]. It is known that the following eigenvalue problem
$$
Au=-\lambda u,\quad u(\partial D)=0
$$
has infinite set of real eigenvalues  $\lambda_{i}\to\infty$ and
$$
0<\lambda_{1}<\lambda_{2}<\cdots<\lambda_{s}<\cdots.
$$
The corresponding eigenfunctions
$$
f_{11},\dots,f_{1n_{1}},\cdots,f_{s1},\dots,f_{sn_{s}},\cdots
$$
form the complete system of functions both in $H_{A}$ and
$L_{2}^{0}(Q):= \{u: u\in L_{2}(Q)\cap u(\partial Q)=0\}$. Here the
number $n_{k}$ is equal to multiplicity of eigenvalue $\lambda_{k}$.

It is often convenient to present the system of eigenfunctions by
one index: $\{f_{n}(z)\}$. The corresponding system of eigenvalues
$\{\lambda_{n}\}$ will be with recurrences. We shall use it too.

We introduce the spectral function

$$
e(x,y,\lambda)=\sum\limits_{\lambda_{j}\leq
\lambda}f_{j}(x)f_{j}(y).
$$

We shall need in the following theorem from the monograph [12].

\proclaim {Theorem 1 ([12].Th.17.5.3)}. There exists such constant
$C_{\alpha}$ that

$$
\sup\limits_{x,y\in D}\sqrt{|D^{\alpha}_{x,y}e(x,y,\lambda)|}\leq
C_{\alpha}\lambda^{(n+|\alpha|)/2}
$$
\endproclaim

Here $\alpha$ is multi-index.
\bigskip
\proclaim
 {\bf Theorem 2}. {\it If the surface $\partial D$ is Lyapunov
surface  and

$$\sup\limits_{(t,z)\in [0,T]\times D, 1\leq i,j\leq d}
\max\left\{|\frac{\partial a_{i}(t,z)}{\partial z_{j}}|,
|\frac{\partial b_{i}(z)}{\partial z_{j}}|,|\frac{\partial
a_{i}(T-t,z)}{\partial t}|\right\}<\infty
$$

then the following relation takes place at $\epsilon\to 0$

$$
P(\epsilon,
z\epsilon)=\exp\left\{-\lambda_{1}\frac{T}{\epsilon^{2}}+\int_{0}^{T}\mu(t)dt\right\}
\sum\limits_{m=1}^{n_{1}}c_{1m}f_{1m}(z)
\left(1+O(\epsilon)\right),\quad \hbox{at}\quad  z\in D,
$$

where
$$
\mu(t)=\sum\limits_{i,j}\left(\frac{1}{2}\sigma_{ij}(0)a_{i}(t,0)a_{j}(t,0)-
\delta_{ij}a_{i}(t,0)a_{j}(t,0)\right)
$$

and $c_{1m}=\int\limits_{D}f_{1m}(z)dz$.}
\endproclaim

\demo {\bf Proof} Make the change of variables and function

$$
x_{i}=z_{i}\epsilon ,\quad u_{1}^{\epsilon}= u_{0}^{\epsilon}
\exp\left\{\epsilon\sum\limits_{k=1}^{d}a_{k}(T-t,0)z_{k}\right\},
$$

 Now we obtain the following parabolic problem for function
$u_{1}^{\epsilon}$

$$
 \frac{\partial u_{1}^{\epsilon}(t,z)}{\partial t}=
\frac{1}{2\epsilon^{2}}\sum\limits_{i,j=1}^{d}\sigma_{ij}(\epsilon
z)\frac{\partial^{2} u_{1}^{\epsilon}(t,z)} {\partial z_{i}\partial
z_{j}}+ \frac{1}{\epsilon}\sum\limits_{i,j}\left(a_{i}(T-t,\epsilon
z) -\frac{1}{2}\sigma_{ij}(\epsilon
z)a_{j}(T-t,0)\right)\frac{\partial u_{1}^{\epsilon}(t,z)}{\partial
z_{i}}
$$
$$
+\sum\limits_{i,j}\left(\frac{1}{2}\sigma_{ij}(\epsilon
z)a_{i}(T-t,0)a_{j}(T-t,0)-
\delta_{ij}a_{i}(T-t,0)a_{j}(T-t,\epsilon z)-\epsilon \frac{\partial
a_{i}(T-t,0)}{\partial t}z_{i}\right)u_{1}^{\epsilon},
 \quad z\in  D;
$$

$$
u_{1}^{\epsilon}(t,z)|_{t=0}=\exp\left\{\epsilon
\sum\limits_{k=1}^{d}a_{k}(T,0)z_{k}\right\};\quad z\in D;\quad
u_{1}^{\epsilon}(t,z)=0\quad z\in
\partial D,\quad 0\leq t\leq T.\eqno (4)
$$

We will construct the  asymptotic expansion of solution for this
initial - boundary problem in the following form

$$
u_{1}^{\epsilon}(t,z)=\sum\limits_{k\geq
0}v_{k}(t,z)\epsilon^{k}.\eqno(5)
$$

Note that  the famous expansion

$$
\exp\left\{\epsilon \sum\limits_{k=1}^{d}a_{k}(T,0)z_{k}\right\}=
1+\epsilon
\sum\limits_{k=1}^{d}a_{k}(T,0)z_{k}+\frac{1}{2!}\left(\epsilon
\sum\limits_{k=1}^{d}a_{k}(T,0)z_{k}\right)^{2}+\cdots,
$$

defines the initial conditions for $v_{k},\quad k\geq 0$:

$$
v_{0}(0,z)=1, \quad v_{1}(0,z)=
\sum\limits_{k=1}^{d}a_{k}(T,0)z_{k},\quad
v_{2}(0,z)=\frac{1}{2}\left(
\sum\limits_{k=1}^{d}a_{k}(T,0)z_{k}\right)^{2}\cdots.
$$

Using the first fragment of Taylor series in zero point under
conditions of theorem we can obtain the following representations

$$
\sigma_{ij}(\epsilon
z)=\sigma_{ij}(0)+\epsilon\sigma^{\epsilon}_{ij}(z),\quad
a_{i}(T-t,\epsilon z)=a_{i}(T-t,0)+\epsilon a_{i}^{\epsilon}(T-t,
z), 1\leq  i,j\leq d\eqno(6)
$$

where

$$
\sup_{z\in D, \epsilon\in [0,1],1\leq i,j\leq
d}|\sigma^{\epsilon}_{ij}(z)|<\infty,\quad \sup_{z\in D, t\in
[0,T],\epsilon\in [0,1],1\leq i\leq d}|a_{i}^{\epsilon}(T-t,
z)|<\infty
$$

 Now, after substitution of (5),(6) to (4) we conclude that the $v_{0}$ satisfies the problem

$$
\frac{\partial v_{0}}{\partial t}=
\frac{1}{2\epsilon^{2}}\left(\sum\limits_{i,j=1}^{d}\sigma_{ij}(0)\frac{\partial^{2}
}{\partial z_{i}\partial z_{j}}\right)v_{0}+\mu(t)v_{0}\eqno(7)
$$
$$
v_{0}|_{\partial D}=0;\quad v_{0}(0,z)=1,\quad z\in D.
$$

Here
$$
\mu(t)=\sum\limits_{i,j}\left(\frac{1}{2}\sigma_{ij}(0)a_{i}(T-t,0)a_{j}(T-t,0)-
\delta_{ij}a_{i}(T-t,0)a_{j}(T-t,0)\right).
$$

Further, let us denote by $B_{\epsilon}(t,z)$ the operator
$C^{2}(D)\to C(D)$, for $f\in C^{2}(D)$ it's defined as follows:

$$
B^{\epsilon}(t,z)f=
$$
$$
\frac{1}{2\epsilon
}\sum\limits_{i,j=1}^{d}\sigma^{\epsilon}_{ij}(z)\frac{\partial^{2}f
}{\partial z_{i}\partial z_{j}}+
\frac{1}{\epsilon}\sum\limits_{i,j}\left(a_{i}(T-t,\epsilon
z)-\frac{1}{2}\sigma_{ij}(\epsilon z)a_{j}(T-t,0)\right)
\frac{\partial f}{\partial z_{i}}+
$$
$$
+\epsilon\sum\limits_{i,j=1}^{d}\left(\frac{1}{2}\sigma^{\epsilon}_{ij}(z)
a_{i}(T-t,0)a_{j}(T-t,0)-\delta_{ij} a_{i}(T-t,0)
a_{j}^{\epsilon}(T-t, z)-\frac{\partial a_{i}(T-t,0)}{\partial
t}z_{i}\right)f=
$$
$$
=:\frac{1}{2\epsilon
}\sum\limits_{i,j=1}^{d}\sigma^{\epsilon}_{ij}(z)\frac{\partial^{2}f
}{\partial z_{i}\partial
z_{j}}+\frac{1}{\epsilon}A^{\epsilon}_{1}(t,z)f+\epsilon
A_{2}^{\epsilon}(t,z).
$$

Now, formally the functions $v_{k}, k\geq 1$ are defined by the
following recurrence system problems
$$
\frac{\partial v_{k}}{\partial t}=
\frac{1}{2\epsilon^{2}}\left(\sum\limits_{i,j=1}^{d}\sigma_{ij}(0)\frac{\partial^{2}
}{\partial z_{i}\partial z_{j}}\right)v_{k}+
B_{\epsilon}(t,z)v_{k-1}\eqno(8)
$$
$$
v_{0}|_{\partial D}=0;\quad
 v_{k}(0,z)=\frac{1}{k!}\left(
\sum\limits_{k=1}^{d}a_{k}(T-t,0)z_{k}\right)^{k}, z\in D.
$$

We shall solve the problems of (7),(8) by method of separation of
variables. According to this method the solutions are defined in the
form

$$
v_{k}(t,z)=\sum\limits_{n\geq 1}q_{k,n}(t)f_{n}(z). \eqno(9)
$$

For definition  of principal number it suffices to construct of the
$v_{0}$. If we substitute (9) at $k=0$ to (7) then we obtain

$$
\sum\limits_{n\geq 1}\left\{-\dot q_{0,n}(t) -
\frac{\lambda_{n}}{\epsilon^{2}}q_{0,n}(t)+\mu(t)q_{0,n}(t)\right\}f_{n}(z)=0.
$$

Set $c_{0,n}=\int_{D}f_{n}(z)dz$ (coefficients of expansion of
indicator of set $D$). The initial condition of $v_{0}$ has the
following stating

$$
v_{0}(0,z)=\sum\limits_{n\geq
1}q_{0,n}(0)f_{n}(z)=\sum\limits_{n\geq
1}c_{0,n}f_{n}(z)=\sum\limits_{l\geq
1}\sum\limits_{m=1}^{n_{l}}c_{0,lm}f_{lm}(z), \quad z\in D.
$$

By definition of system of functions $\{f_{n}(z)\}$, now we have
 the system of ordinary differential equations

$$
\dot
q_{0,n}(t)+\left(\frac{\lambda_{n}}{\epsilon^{2}}-\mu(t)\right)q_{0,n}(t)=0,\quad
q_{0,n}(0)=c_{0,n}.
$$

From the latter one we have

$$
q_{0,n}(t)=c_{0,n}\exp\left\{-\frac{\lambda_{n}}{\epsilon^{2}}t+
\int_{0}^{t}\mu(s)ds\right\}.
$$

Set
$$
A_{0}=\sup\limits_{\epsilon\leq 1,z\in D;
i,j}|\sigma_{ij}^{\epsilon}(z)|,\quad L_{0}=\sum\limits_{l\geq
1,1\leq m\leq n_{l}}\left(c_{0,ml}\right)^{2}.
$$
$$
A_{1}=\sup\limits_{0\leq\epsilon\leq 1, z\in D,t\in [0,T];i,j}
\left|a_{i}(T-t,\epsilon z)-\frac{1}{2}\sigma_{ij}(\epsilon
z)a_{j}(T-t,0)\right|.
$$
$$
A_{2}=
$$
$$
=\sup\limits_{0\leq\epsilon\leq 1, z\in D,t\in [0,T]; i,j}
\left|\frac{1}{2}\sigma^{\epsilon}_{ij}(z)
a_{i}(T-t,0)a_{j}(T-t,0)-\delta_{ij} a_{i}(T-t,0)
a_{j}^{\epsilon}(T-t, z)-\frac{\partial a_{i}(T-t,0)}{\partial
t}z_{i}\right|
$$

We have the following relations for eigenvalues $\lambda_{l}$

$$
k_{1}l^{2/d}\leq \lambda_{l}\leq k_{2}l^{2/d},\quad
\max(k_{1},k_{2})<\infty
$$,

Applying Cauchy-Bunyakovskii inequality, Theorem 1 and the latter
one, we get

$$
\left|\sum\limits_{i,j}
a^{\epsilon}_{i,j}(z)\frac{\partial^{2}v_{0}}{\partial z_{i}\partial
z_{j}}\right|=\left|\sum\limits_{l}\exp\left( -
\lambda_{l}t\epsilon^{-2}+
\int_{0}^{t}\mu(s)ds\right)\sum\limits_{m=1}^{n_{l}}c_{0,ml}\sum\limits_{i,j}a^{\epsilon}_{i,j}(z)
\frac{\partial^{2}f_{ml}(z)} {\partial z_{i}\partial
z_{j}}\right|\leq
$$
$$
\leq A_{0} d\sum\limits_{l}\exp\left( - \lambda_{l}t\epsilon^{-2}+
\int_{0}^{t}\mu(s)ds\right)\left(\sum\limits_{m=1}^{n_{l}}\left(c_{0,ml}\right)^{2}\right)
^{\frac{1}{2}}\left(\sum\limits_{m=1}^{n_{l}}\sum\limits_{i,j}\left(\frac{\partial^{2}f_{ml}(z)}
{\partial z_{i}\partial z_{j}}\right)^{2}\right)^{\frac{1}{2}}\leq
$$
$$
\leq A_{0} d C_{2,2} L_{0} \sum\limits_{l\geq
1}\exp\left(-\frac{\lambda_{l}t}{\epsilon^{2}}+\int\limits_{0}^{t}\mu(s)ds\right)\lambda_{l}^{\frac{d}{2}+2}
\leq
\exp\left(-\frac{\lambda_{1}t}{\epsilon^{2}}\right)K_{0}.\eqno(10)
$$

Here $K_{0}<\infty$.

Reasoning similarly we convince ourselves that for other parts of
$B^{\epsilon}(t,z)v_{0}$ the following estimations take place

$$
|A_{1}^{\epsilon}(t,z)v_{0}|\leq A_{1}d
C_{1,1}L_{0}\sum\limits_{l\geq
1}\exp\left(-\frac{\lambda_{l}t}{\epsilon^{2}}+\int\limits_{0}^{t}\mu(s)ds\right)\lambda_{l}^{\frac{d}{2}+1}
\leq
\exp\left(-\frac{\lambda_{1}t}{\epsilon^{2}}\right)K_{0,1};\eqno(11)
$$
$$
|A_{2}^{\epsilon}(t,z)v_{0}|\leq A_{2} d
C_{0,0}L_{0}\sum\limits_{l\geq
1}\exp\left(-\frac{\lambda_{l}t}{\epsilon^{2}}+\int\limits_{0}^{t}\mu(s)ds\right)\lambda_{l}^{\frac{d}{2}}
\leq
\exp\left(-\frac{\lambda_{1}t}{\epsilon^{2}}\right)K_{0,2},\eqno(12)
$$

where $\max\{K_{0,1},K_{0,2}\}<\infty.$

Now let us estimate the coefficients  $\beta^{\epsilon}_{n}(t)$ of
expansion of $B^{\epsilon}(t,z)v_{0}$ by system $\{f_{n}\}_{n\geq
1}$. Applying (10)-(12) and Cauchy-Bunyakovskii inequality, we get

$$
|\beta^{\epsilon}_{n}(t)|=|\int\limits_{D}B^{\epsilon}(t,z)v_{0}(t,z)f_{n}(z)dz|\leq
$$
$$
\leq
\left(\int\limits_{D}\left(B^{\epsilon}(t,z)v_{0}\right)^{2}dz\right)^{\frac{1}{2}}
\left(\int\limits_{D}f_{n}^{2}(z)dz\right)^{\frac{1}{2}}\leq
$$
$$
\leq
\exp(-\lambda_{1}t\epsilon^{-2})\left(\frac{K_{0}+K_{0,1}}{\epsilon}+\epsilon
K_{0,2}\right)|D|.
$$

The latter one now gives

$$
|\int\limits_{0}^{t}\beta^{\epsilon}_{n}(s)ds|\leq \epsilon
\gamma_{\epsilon}(t),\eqno(13)
$$

where
$$
\sup\limits_{0 \leq \epsilon \leq 1, t\in
[0,T]}\gamma_{\epsilon}(t)<\infty.
$$

Finally, let us estimate the difference
$r^{\epsilon}(t,z)=u^{\epsilon}_{1}(t,z)-v_{0}(t,z)$. By definition,
$r^{\epsilon}(t,z)$ is solution of the following problem

$$
\frac{\partial r^{\epsilon}}{\partial t}=
\frac{1}{2\epsilon^{2}}\left(\sum\limits_{i,j=1}^{d}\sigma_{ij}(0)\frac{\partial^{2}
}{\partial z_{i}\partial z_{j}}\right)r^{\epsilon}+
B_{\epsilon}(T-t,z)v_{0}\quad z\in  D;\eqno(14)
$$

$$
r^{\epsilon}(t,z)|_{t=0}=\exp\left\{\epsilon
\sum\limits_{k=1}^{d}a_{k}(T,0)z_{k}\right\}-1;\quad z\in D;\quad
r^{\epsilon}(t,z)=0\quad z\in
\partial D,\quad 0\leq t\leq T.
$$

It is clear that $r^{\epsilon}(0,z)$ we can present as $\epsilon
r_{1}^{\epsilon}(0,z)$, where $r_{1}^{\epsilon}(0,z)$ is uniform
bounded function of variables $\epsilon\in [0,1]$ and $z\in D$.  So,
the coefficients of expansion this function by system $\{f_{n}(z)\}$
have the following forms

$$
\int\limits_{D}r^{\epsilon}(0,z)f_{n}(z)dz= \epsilon
\mu_{n}^{\epsilon},\quad  \hbox{where}\quad \sup\limits_{0\leq
\epsilon\leq 1}\sum\limits_{n\geq 1}
\left(\mu_{n}^{\epsilon}\right)^{2}=M<\infty.\eqno(13)
$$

Now we have the solution of (14) in the following form

$$
r^{\epsilon}(t,z)=\epsilon\sum\limits_{n\geq
1}\mu_{n}^{\epsilon}\exp\{-\lambda_{n}t\epsilon^{-2}+\int\limits\beta_{n}^{\epsilon}(s)ds\}f_{n}(z)
$$

Applying latter one ,(13),(15), Theorem 1 and Cauchy-Bunyakovskii
inequality we get at $t>0$

$$
|r^{\epsilon}(t,z)\epsilon ^{-1}|\leq \left(\sum\limits_{n\geq
1}(\mu_{n}^{\epsilon})^{2}\right)^{\frac{1}{2}}C_{0,0}\sum\limits_{n\geq
1}\exp\{-\lambda_{n}t\epsilon^{-2}+\int\limits\beta_{n}^{\epsilon}(s)ds\}\lambda_{n}^{\frac{d}{2}}\}\leq
$$
$$
\leq M C_{0,0}\exp\{-\lambda_{1}t\epsilon^{-2}\} K_{0,3},\quad
\hbox{where}\quad K_{0,3}<\infty.
$$

The proof of theorem is completed.

\enddemo

\proclaim
 {\bf Remark 1}According to the above system of problems for definition of the functions $v_{k}, k\geq 1$,
 we outline the construction of coefficients $q_{k,n})(t)$ for the series (8):

$$
\dot q_{k,n}(t)= +\left(\frac{\lambda_{n}}{\epsilon^{2}}+
\mu^{\epsilon}_{k-1,n}(t)\right)q_{k,n}(t),
$$
$$
q_{k,n}(0)=\int_{D}v_{k}(0,z)f_{n}(z)dz= \frac{1}{k!}\int_{D}\left(
\sum\limits_{m=1}^{d}a_{m}(T,0)z_{m}\right)^{m}f_{n}(z)dz
$$

Here $\mu^{\epsilon}_{k-1,n}(t)= \int\limits_{D}f_{n}(z)
B^{\epsilon}(t,z)v_{k-1}(t,z)dz.$
\endproclaim

\proclaim {\bf Remark 2} Theorem 2 is coordinated with results of
works [6-8] where the principal member of small deviations in ball
are investigated for more simple SDE.
\endproclaim

\bigskip

{\bf II.} {\bf The rarefaction of set of diffusion processes by
boundaries of small domains.}
\bigskip

The following problem was investigated in works[13,14]. Let a set
identical diffusion random processes start at the initial time from
the different points of domain $D$. These processes are diffusion
processes with absorbtion on the boundary $\partial D$. We are
interested in distribution of the number yet absorbed at the moment
$T$. The initial number and initial position of diffusion processes
are defined either a random Poisson measure[14] or deterministic
measure [13]. The proved limits theorems  described the situation
when $T\to \infty$ and initial number of diffusion processes
depended on $T$ and it increased at the rise of $T$. The  role of
normalizing function played principal member of asymptote of
solution of according parabolic problem at  $T\to \infty$.

Henceforth  we shall assume that considered diffusion processes
satisfy of the  SDE (2) with different initial points.

Now we consider the situation when initial number of absorbing
diffusion processes in small domain $\epsilon D$ depends on
$\epsilon\to 0$ and it increase under the condition of decrease of
$\epsilon$.
 It is not hard to show, that now normalizing function is the
 principal member of parabolic problem (3) at $\epsilon\to 0$.

The proofs of stated below theorems repeat the proofs of according
theorems from [13,14] almost word for word.

 We will denote by ~$\eta(\epsilon,T)$~ the number of remaining
processes in the region ~$\epsilon D$~ at the moment ~$T$.

We will also assume that $\sigma$-additive  measure $\nu$ is given
on the $\Sigma_{\nu}$- algebra sets from $D,\quad \nu(D)<\infty.$
All eigenfunctions $f_{ij}:D\to R^{1}$ are
$(\Sigma_{\nu},\Sigma_{Y})$ measurable. Here $\Sigma_{Y}$ is system
of Borel sets from $R^{1}$. Let $\Rightarrow$ denote the  weak
convergence of random values or measures.

 At the beginning we assume that initial number and position of
 diffusion processes are defined by deterministic measure
 $N(\epsilon B,\epsilon), B\in D$. Thus, $N(\epsilon B,\epsilon)$ is
 equal to number of starting points in the set $\epsilon B$.

Let us denote by ~$\nu_{\epsilon}(\cdot)$ the measure

$$
\nu_{\epsilon}(\epsilon
B)=\exp\left(-\frac{T\lambda_{1}}{2\epsilon^{2}}\right)N(\epsilon
B,\epsilon).
$$

where $B\in\Sigma_{\nu}$ .

By definition of measure ~$\nu_{\epsilon}(\cdot)$, we have
$$
d\nu_{\epsilon}(x)= \cases
\exp\left(-\frac{T\lambda_{1}}{2\epsilon^{2}}\right)
,&  \hbox{if}\quad x=x_{k},\quad k=1,\cdots, N(\epsilon D,\epsilon) \\
0 ,& \hbox{otherwise}.
\endcases
$$

\bigskip
\proclaim {\bf Theorem 3} {\it Under the assumptions of the Theorem
2 let the  $N(\epsilon \cdot, \epsilon)$ satisfies the condition

$$
\nu_{\epsilon}(\epsilon~
\cdot)\mathop{\Rightarrow}\limits_{\epsilon\to 0}\nu(\cdot).
$$

Then $\eta(\epsilon,T)\Rightarrow\eta(T)$ if $\epsilon \to 0$ where
$\eta(T)$ has Poisson distribution function with parameter

$$a(T)=\exp\left(\int\limits_{0}^{T}\mu(s)ds\right)\int\limits_{D}
F(z)d\nu(z),
$$

where $F(z)=\sum\limits_{i=1}^{n_{1}}f_{1i}(z)c_{1i},\quad
c_{1i}=\int\limits_{Q}f_{1i}(z)d z $

and $\mu(t)$ is the function from Theorem 2.}
\endproclaim

Now we consider the case when the initial number and positions of
processes are defined by the random Poisson measure
$\mu(\cdot,\epsilon)$ in $\epsilon D$:
$$
P(\mu(\epsilon A,\epsilon)=k)=\frac{m^{k}(\epsilon
A,\epsilon)}{k!}e^{-m(\epsilon A,\epsilon)},
$$
where~$m(\epsilon ~\cdot,\epsilon)$ is finitely additive positive
measure on ~$\epsilon D$ for fixed $\epsilon$.

We assign

$$
g(\epsilon)=\exp\left(-\frac{T\lambda_{1}}{2\epsilon^{2}}\right).
$$

\bigskip
\proclaim
 {\bf Theorem 4}{\it Under the assumptions of the Theorem 2
we suppose that $m(\epsilon \cdot, \epsilon)$ holds the condition

$$
\lim\limits_{\epsilon\to 0} m(\epsilon B,\epsilon)g(\epsilon)=
\nu(B),\quad B\in \Sigma_{\nu}.
$$

Then $\eta(\epsilon,T)\Rightarrow\eta(T)$ if $\epsilon\to 0$ where
$\eta(T)$ has the Poisson distribution function with the parameter
$a(T)$ from Theorem 3.}
\endproclaim

\bigskip

 {\bf References}

\bigskip
1. Graham R.,Path integral  formulation of general diffusion
processes, Z.Phys.(1979),B 26,pp.281-290.

2. Onsager L. and Machlup S. Fluctuation and irreversible processes,
I,II, Phys.Rev.(1953) 91,pp.1505-1512,1512-1515.

3. Li W. V.,Shao Q.-M., Gaussian processes:inequalities, small ball
probabilities and applications, in : Stochastic Processes:Theory and
Methods, in : Handbook of Statistics, vol.19, 2001, pp. 533-597.

4. Lifshits M.A., Asymptotic behavior of small ball probabilities,
in Probab. Theory and Math.Statist., Proc. VII International Vilnius
Conference (1998), pp. 453-468.

5. Lifshits M., Simon T., Small deviations for fractional stable
processes, Ann.  I. H. Poincare - PR 41 (2005) pp. 725-752.

6. Mogulskii A.A, The method of Fourier for determination of
asymptotics of small deviations of Wiener process, Siberian Math.
Journ. (1982),v.22,no.3,pp.161-174.

7. Fujita T. and Kotani S., The Onsager - Machlup Function for
diffusion processes, J.Math.Kyoto Uneversity.-
1982.-vol.22,no.22.pp.131-153.

8. Zeitoni O., On the Onsager-Machlup functional of diffusion
processes around non $C^{2}$ curves, Ann. Probab.(1989),vol.17,
no.3, pp.1037-1054.

9. Gasanenko V.A., The total asymptotic expansion of sojourn
probability of diffusion process in thin domain with moving
boundaries, Ukraine Math. Journ. (1999),v.51, no. 9, pp.1155-1164.

10. Gasanenko V.A., The jump like processes in thin domain, Analytic
questions of stochastic system, Kyiv:Institute of Mathematics
(1992), pp. 4-9.

11. Mihlin S.G. Partial differential linear  equations (1977),
Vyshaij shkola, Moskow, 431.

12.L.H$\ddot o$rmander, The analysis of Linear Partial differential
Operators III (1985), Spinger-Verlag.

13.Fedullo A., Gasanenko V.A., Limit theorems for rarefaction of set
of diffusion processes by boundaries, Theory of Stochastic Processes
vol. 11(27), no.1-2,2005, pp.23-29.

14.Fedullo A., Gasanenko V.A.,Limit theorems for number of diffusion
processes, which did not absorb by boundaries, Central European
Journal of Mathematics 4(4), 2006, pp.624-634.
\enddocument